\documentclass[11pt,oneside,a4paper]{amsart}
\usepackage{graphicx}
\usepackage{subfigure}
\usepackage{amsmath,amssymb,amsthm,bbm,marvosym}
\usepackage{enumitem}
\usepackage[]{units}
\setenumerate{leftmargin=*}

\usepackage{fullpage}

\let\subset\subseteq 
\let\eps\varepsilon
\let\rho\varrho

\def\NN{\mathbb{N}}

\newcommand{\JUSTIFY}[1]{\fbox{\tiny{#1}}\quad}

\def\differential{\mathsf{d}} 
\def\Pr{\mathbf{P}} 
\def\Exp{\mathbf{E}} 
\def\Poi{\mathsf{Poi}} 
\def\Uni{\mathsf{Uni}} 
\def\Dirac{\mathsf{Dirac}} 
\def\RG{\mathbb{G}} 
\def\G{\mathbf{H}} 

\def\cD{\mathcal{D}}
\def\fX{\mathfrak{X}}
\def\mindeg{\mathrm{mindeg}}
\def\cutdist{\delta_{\square}}
\def\cutnormdist{d_{\square}}
\def\previouslyweight{length{}}
\def\previouslyweights{lengths{}}

\newtheorem{theorem}{Theorem}
\newtheorem{lemma}[theorem] {Lemma}   
   
\newtheorem{conjecture}[theorem] {Conjecture}  
  
\newtheorem{fact}[theorem]{Fact}

\newtheorem{claimNewNumbering}{Claim}
\newtheorem{claim}[claimNewNumbering]{Claim}
\newtheorem*{claim*}{Claim}
\newtheorem*{fact*}{Fact} 
\newtheorem{remark}[theorem] {Remark} 
\newtheorem{mainsetup}[theorem] {Setup} 
\newtheorem{definition}[theorem] {Definition} 
\newtheorem*{remark*} {Remark}

\newcommand{\By}[2]{\overset{\mbox{\tiny{#1}}}{#2}}
\newcommand{\ByRef}[2]{   \By{\eqref{#1}}{#2} }

\newcommand{\leBy}[1]{    \By{#1}{\le} }
\newcommand{\geBy}[1]{    \By{#1}{\ge} }

\newcommand{\leByRef}[1]{ \ByRef{#1}{\le} }

\renewcommand{\leq}{\leqslant}
\renewcommand{\le}{\leqslant}
\renewcommand{\geq}{\geqslant}
\renewcommand{\ge}{\geqslant}

\renewcommand{\epsilon}{\varepsilon}

\title{Random minimum spanning tree and dense graph limits}
\author{Jan Hladký}
\address{Institute of Computer Science of the Czech Academy of Sciences, Pod Vod\'{a}renskou v\v{e}\v{z}\'{\i} 2, 182~00 Prague, Czechia. With institutional support RVO:67985807.}
\email{hladky@cs.cas.cz}
\author{Gopal Viswanathan}
\address{Charles University, Faculty of Mathematics and Physics, Czechia.}
\thanks{Research supported by Czech Science Foundation Project 21-21762X}
\subjclass[2020]{05C80}
\keywords{random minimum spanning tree, dense graph limits, graphons}
\date{}

\begin{document}
\begin{abstract}
	A theorem of Frieze from~1985 asserts that the total \previouslyweight\ of the minimum spanning tree of the complete graph $K_n$ whose edges get independent \previouslyweights\ from the distribution $\Uni[0,1]$ converges to Apéry's constant in probability, as $n\to\infty$. We generalize this result to sequences of graphs $G_n$ that converge to a graphon $W$. Further, we allow the \previouslyweights\  of the edges to be drawn from different distributions (subject to moderate conditions). The limiting total \previouslyweight\ $\kappa(W)$ of the minimum spanning tree is expressed in terms of a certain branching process defined on $W$, which was studied previously by Bollobás, Janson and Riordan in connection with the giant component in inhomogeneous random graphs.
\end{abstract}
\maketitle

\section{Introduction}
The main contribution of the paper is a continuity result of the random minimum spanning tree problem with respect to the topology of dense graph limits. Let us first introduce basic graph-theoretic notation and then summarize the state of the art in the area of random minimum spanning trees. If $G$ is a graph (finite, undirected, simple), we write $V(G)$ and $E(G)$ for its vertex set and its edge set, respectively. Further, we write $v(G)=|V(G)|$ and $e(G)=|E(G)|$. If $G$ is connected and $w:E(G)\rightarrow [0,\infty)$ is a function (which we call the \emph{\previouslyweight\ 
 function}), then the \emph{\previouslyweight}\footnote{\label{foot:weightlenght}While this is commonly termed the `weight' of an edge or a spanning tree, we opt for `length' to mitigate potential ambiguity. The word `weight' will also be used in the context of the numerical values associated with edges in certain weighted graphs, thus necessitating a distinct term.} of a spanning tree $T$ is 
\begin{equation}\label{eq:MSTdef}
w(T):=\sum_{e\in E(T)}w(e)\;.
\end{equation}
A \emph{minimum spanning tree} is a spanning tree of $G$ which minimizes $w(T)$. We write $MST(G):=w(T)$, where $T$ is a minimum spanning tree of $G$. A minimum spanning tree has many nice features which also allow very simple and efficient algorithms which in turn can help analyze mathematical problems concerning a minimum spanning tree. These are Borůvka's algorithm from 1926, Jarník's algorithm from 1930 (rediscovered and popularized by Prim in 1957 and by Dijkstra in 1959) and two algorithms which appeared in a 1956 paper by Kruskal \cite{KruskalAlgorithms}. One of those two algorithms turns out to be much more useful than the other and has been called `Kruskal's algorithm'. We shall recall (and then use) this algorithm below.

We now move to the context of a random minimum spanning tree. Here, the \previouslyweight\ of each edge $e$ of the graph $G$ is chosen independently from a certain probability distribution $D_e$ on $[0,\infty)$. Then the \previouslyweight\ of the minimum spanning tree is a random variable, whose distribution we denote by $MST(G,(D_e)_{e\in E(G)})$. If all the distributions are equal to the same distribution $D$, we simply write $MST(G,D)$. The first limit theorem about the random minimum spanning tree was obtained by Frieze~\cite{MR770868} for the complete graph $K_n$ and the uniform distribution on $[0,1]$ on each edge. We state  a slight strengthening obtained shortly after by Steele~\cite{MR0905183}.
\begin{theorem}\label{thm:frieze}
The sequence $MST(K_n,\Uni[0,1])$ converges in probability to Apéry's constant $\zeta(3)\approx 1.202$ as $n\to \infty$.
\end{theorem}
It is worth noting that Theorem~\ref{thm:frieze} does not involve any rescaling. That means that as $n$ grows, the increase in the number of summands in~\eqref{eq:MSTdef} is compensated by the fact that these summands are typically becoming smaller (thanks to a bigger pool of $\Uni[0,1]$ random variables to choose from).

In fact, Frieze's result is more general. Suppose that $D$ is a probability distribution on $[0,\infty)$. Let $F$ be the cumulative distribution function of $D$ and suppose that $F$ is differentiable at~0. We write $F'(0)$ for the derivative and assume that $F'(0)>0$.
\begin{theorem}\label{thm:frieze2}
With the assumption above, the sequence $MST(K_n,D)$ converges in probability  to $\zeta(3)/F'(0)$ as $n\to \infty$.
\end{theorem}
Frieze's result has been extended in a number of ways. Here we recall a particularly relevant strengthening by Frieze and McDiarmid~\cite{MR1054012}. Suppose that $H$ is a graph, $(D_e)_{e\in E(H)}$ are probability distributions on $[0,\infty)$ with cumulative distribution functions $(F_e)_{e\in E(H)}$ whose derivatives at~0 are $(F'_e(0))_{e\in E(H)}$, and that they are all positive. We allow $H$ to have self-loops.  We say that $(H,(D_e)_{e\in E(H)})$ is a \emph{$\Delta$-regular template} (for some $\Delta\in(0,\infty)$) if for every $v\in V(H)$ we have $\sum_{e\in E(H):e\ni v}F'_e(0)=\Delta$. Here, a self-loop at $v$ contributes to the sum only once.

Let $\mathfrak{H}=(H,(D_e)_{e\in E(H)})$ be a graph whose edges are equipped with probability distributions. Let $n\in\NN$. Then the \emph{$n$-blow-up of $\mathfrak{H}$} is a graph $G$ whose edges are equipped with probability distributions, defined as follows. The vertex set of $G$ is $[n]\times V(H)$. A pair of distinct vertices $(i,u)$ and $(j,v)$ forms an edge of $G$ if and only if $uv\in E(H)$. In such a case the probability distribution on that edge is $D_{uv}$.

The result of Frieze and McDiarmid reads as follows.
\begin{theorem}\label{thm:FriezeMcDiarmid}
Suppose that $\mathfrak{H}=(H,(D_e)_{e\in E(H)})$ is a $\Delta$-regular template in which $H$ is a connected graph. Let $(G_n,\cD_n)$ be a sequence of $n$-blow-ups of $\mathfrak{H}$. Then $MST(G_n,\cD_n)$ converges in probability to $\zeta(3)/\Delta$.
\end{theorem}
Theorem~\ref{thm:FriezeMcDiarmid} indeed generalizes Theorem~\ref{thm:frieze2}. To this end, it is enough to observe that $K_n$ is the $n$-blow-up of a single vertex with a self-loop.
Note that other results about the random minimum spanning tree, as far as we could find, also require regularity of the degrees (see for example~\cite{MR1721947}).

\subsection{Our result}
Our main result, Theorem~\ref{thm:mainfull} below, is an extension of Frieze's theorem to sequences of dense graphs whose edges are equipped with general distributions on their edges. To this end, we use the theory of dense graph limits. This theory, initiated in~\cite{Lovasz2006,MR2455626}, compactifies the space of finite graphs by objects called `graphons' (and in our case also by somewhat more general `kernels') which are certain analytic counterparts to graphs. This analytic view offers powerful additional tools that have led to number of breakthroughs in extremal and random graph theory. We use only basic aspects of the theory, which we summarize in Section~\ref{sss:Graphons}. More specifically, Theorem~\ref{thm:mainfull} asserts that if we weight the edges by the derivative of the cumulative distribution functions at~0 and this sequence of weighted graphs converges in the cut distance to a kernel $W$, then the \previouslyweights\ of the random minimum spanning trees converge in probability to a certain constant $\kappa(W) \in (0,\infty)$. As we show in Section~\ref{ssec:MAINvsFrMcD}, this result implies Theorem~\ref{thm:FriezeMcDiarmid} but importantly it also applies to sequences of graphs which are not regular. This is the first result, as far as we know, about the random minimum spanning tree in a nonregular setting. Also, the compactness property of the cut distance topology (stated in Theorem~\ref{thm:compactness} below) means that Theorem~\ref{thm:mainfull} provides an asymptotic description of the random minimum spanning tree of \emph{every} sequence of dense graphs with 
edges equipped with independent lengths coming from the possibly different distributions.
Recall that the term `dense graph sequence' traditionally refers to sequences of graphs $G_1,G_2,\ldots$ where $\liminf_n \frac{e(G_n)}{v(G_n)^2}>0$. In our setting we in addition require `robustness', which  prohibits sparse cuts. As we show in Section~\ref{ssec:necessityofassumptions}, this strengthening is necessary.

\subsubsection{Decent families of distributions}\label{sss:decent}
We define the class of probability distributions which we allow to use as edge \previouslyweights. We only work with distributions supported on $[0,+\infty)$. While we allow different distributions for different edges, we require certain uniform differentiability of their cumulative distributions functions at~0. This is summarized in the next definition.
\begin{definition}\label{def:decent}
Suppose that $(\cD_n)_n=(D_{n,i})_{i\in I_n}$ is a sequence of families of probability distributions on $[0,+\infty)$. We say that $(\cD_n)_n$ is \emph{decent} if the following conditions are met.
\begin{enumerate}[label=(d\arabic*)]
\item For each $n\in\NN$ and $i\in I_n$, let $F_{n,i}$ be the cumulative distribution functions of $D_{n,i}$. Suppose that $F_{n,i}$ is differentiable at~0. We write $F_{n,i}'(0)$ for the derivative.
\item\label{en:decent3} For every $\eps>0$, there exists $n_0\in \NN$ and $\delta>0$, such that for each $n\ge n_0$, $i\in I_n$ and each $\alpha\in(0,\delta)$, we have $|F_{n,i}'(0)-\frac{F_{n,i}(\alpha)}{\alpha}|<\eps$.
\end{enumerate}
\end{definition}
In Section~\ref{ssec:necessityofassumptions} we discuss why this is the right definition in our setting. It is routine to deduce from~\ref{en:decent3} that 
\begin{equation}\label{eq:sup}
\sup_{n\in \NN, i\in I_n} F_{n,i}'(0)<\infty\;.
\end{equation}

\subsubsection{Derivatives of the cumulative distribution functions as edge weights}\label{sss:derivativesasedgeweights}
In the simplest instance, our Theorem~\ref{thm:mainfull} below says that if $(H_n)_{n\in \NN}$ is a sequence of well-connected (using a definition of `robustness' below) graphs that converges to a graphon $W$ in the cut distance and we equip the edges of each graph with $\Uni[0,1]$-distributions, then the \previouslyweights\ of the minimum spanning trees converge in probability to a certain constant $\kappa(W)\in (0,\infty)$. Our Theorem~\ref{thm:mainfull}, however, allows different probability distributions to be placed on different edges of $H_n$. As we learned in Theorem~\ref{thm:frieze2} and Theorem~\ref{thm:FriezeMcDiarmid}, the strength of an edge should be the derivative of the cumulative distribution function at~0. This leads to the following definition. Suppose that $H$ is a graph and $\cD=(D_{e})_{e\in E(H)}$ is a family of probability distributions whose cumulative distribution functions $(F_{e})_{e\in E(H)}$ are differentiable at~0. Then the \emph{weighted graph corresponding to $H$ and $\cD$} is the graph $H$ where on each edge $e\in E(H)$ we put weight $F_e'(0)$ (which we assume is positive).

Next, we define an expansion-like notion for weighted graphs, which we call robustness.
If $G$ is a weighted graph whose weight function is $f$, $f:E(G)\to[0,\infty)$. Then for $X,Y\subset V(G)$ we write $e_G(X,Y)=\sum_{x\in X,y\in Y}f(x,y)$.
\begin{definition}
For $\rho>0$, we say that a weighted graph $G$ is \emph{$\rho$-robust}, if for every $U\subset V(G)$ we have $e_G(U,V(G)\setminus U)\ge\rho |U|(v(G)-|U|)$. 
\end{definition}

\subsubsection{Graphons, kernels and cut distance convergence}\label{sss:Graphons}
Throughout the paper, $(\Omega,\mu)$ is a separable atomless probability space with an implicit sigma-algebra. A \emph{kernel} is a symmetric nonnegative function $W\in L^\infty(\Omega^2)$. We say that $W$ is a \emph{graphon} if the range of $W$ is a subset of $[0,1]$. 

Our notation mostly follows~\cite{Lovasz2012}, to which we refer the reader for details. The \emph{minimum degree} of a kernel $W$ is defined as the essential infimum of the degrees of $W$, $\mindeg(W):=\mathrm{essinf}\left\{\int_{y\in\Omega}W(x,y)\differential\nu(y):x\in\Omega\right\}$. For $d\ge 0$, we say that $W$ is \emph{$d$-regular} if we have that $\int_{y\in\Omega}W(x,y)\differential\nu(y)=d$ for almost every $x\in \Omega$. 

Next, we describe graphon representations of finite graphs. Suppose that $G$ is a graph of order $n$. We partition $\Omega$ into sets $\{\Omega_v\}_{v\in V(G)}$ of measure $\frac{1}{n}$ each. We then define the \emph{representation} of $G$ as a graphon $W_G$, which on each square $\Omega_u\times \Omega_v$ is equal to~1 or~0 depending on whether $uv\in E(G)$ or not. Representations of graphs extend to weighted graphs. In that case instead of value~1 for an edge we use the weight of that edges. As a consequence, the representation is not necessarily a $\{0,1\}$-valued graphon but rather a more general step-kernel. 

Note that the above representation is not unique as it depends on the partition $\{\Omega_v\}_{v\in V(G)}$. This is not an issue in our context and reflects the fact that we work with graphs modulo isomorphism. 

The most favorable topology for the theory of dense graph limits is that of the cut distance. Its construction has two steps. Suppose that $U$ and $W$ are two kernels. Then we define the \emph{cut norm distance} between $U$ and $W$ as $\cutnormdist(U,W):=\sup_{S,T\subset \Omega}\left|\int_{S\times T}U-\int_{S\times T}W\right|$. Next, we define the \emph{cut distance} between $U$ and $W$ as $\cutdist(U,W):=\inf_{\pi}\cutnormdist(U,W^\pi)$, where $\pi$ ranges through all measure-preserving bijections on $\Omega$ and $W^\pi$ is a kernel defined by $(W^\pi)(x,y)=W(\pi(x),\pi(y))$. Informally, the role of the measure-preserving bijections $\pi$ is to transfer to the factor-space of isomorphism classes graphs (or graphons and kernels).

Note that for a weighted graph $G$ and a kernel $U$, the distance $\cutdist(W_G,U)$ does not depend on the choice of the $\{\Omega_v\}_{v\in V(G)}$ used to represent $G$. In particular, the cut distance allows us to measure distances between a (weighted) graph and a kernel or between two graphs. 

The key result in the area is that of the compactness of the cut distance topology due to Lovász and Szegedy, \cite{Lovasz2006}. The version stated below for kernels uniformly bounded in the $L^\infty$-norm follows simply by rescaling.
\begin{theorem}\label{thm:compactness}
Let $C>0$ be arbitrary. Then the space of kernels whose $L^\infty$-norm is at most~$C$ equipped with the cut distance is compact.

In particular, if $G_1,G_2,\ldots$ is a sequence of weighted graphs whose weights are bounded by~$C$ then there exists a kernel $W$ and a subsequence $G_{n_1}, G_{n_2},\ldots$ that converges to $W$ in the cut distance.
\end{theorem}

We note that in this paper we could work with the cut distance as a black box. That is, we never actually utilize its precise definition above (with the only exception of one easy step in the proof of Fact~\ref{fact:kappafinite}, which had been known previously). This is because all the favourable properties of the cut distance convergence are contained in Lemma~\ref{lem:FirstLemma} which is derived from the machinery of~\cite{MR2599196}.

The last concept we use from the theory of graph limits is that of homomorphism densities. Suppose that $F$ is a graph (unweighted) and $W$ is a kernel. Then the \emph{homomorphism density of $F$ in $W$} is defined as
\begin{equation}\label{eq:homdens}
t(F,W):=\int \prod_{uv\in E(F)}W(x_u,x_v)
\;
\differential\nu^{\otimes V(F)}\left((x_u)_{u\in V(F)}\right)\;.
\end{equation}

\subsubsection{Statement of the main result}
Theorem~\ref{thm:mainfull} is formulated for sequences of graphs whose edges are equipped with probability distributions and which satisfy additional mild requirements. We fix this setup below.
\begin{mainsetup}\label{mainsetup}
Suppose that $(H_n)_{n\in \NN}$ is a sequence of graphs of growing orders whose edges are equipped with probability distributions $\cD_n=(D_{n,e})_{e\in E(H_n)}$ on $[0,+\infty)$. Suppose that the sequence $(\cD_n)_n$ is decent. For each $n\in\NN$, let $G_n$ be the weighted graph corresponding to $H_n$ and $\cD_n$. Suppose that there exist $\rho>0$ such that all graphs $G_n$ are $\rho$-robust.
\end{mainsetup}

\begin{theorem}\label{thm:mainfull}
Suppose that $(H_n)_{n\in \NN}$, $(\cD_n)_{n\in\NN}=\big((D_{n,e})_{e\in E(H_n)}\big)_{n\in\NN}$, and $(G_n)_{n\in\NN}$ are as in Setup~\ref{mainsetup}, and that the sequence $(G_n)_{n\in\NN}$ converges to a kernel $W$ in the cut distance. Then the sequence $(MST(H_n,\cD_n))_n$ converges in probability to a constant $\kappa(W)\in (0,\infty)$ defined in Section~\ref{ssec:meaningofconstant}.
\end{theorem}

Section~\ref{ssec:meaningofconstant} includes further properties of $\kappa(W)$. In particular, we show that for the constant-1 kernel we indeed have $\kappa(\mathbbm{1})=\zeta(3)$. Furthermore, we introduce the notion `fractional multiple' which allows to relate the value of $\kappa$ in many cases. In particular, we show in Section~\ref{ssec:MAINvsFrMcD} that our result also implies Theorem~\ref{thm:FriezeMcDiarmid}.

\subsubsection{Necessity of the assumptions}\label{ssec:necessityofassumptions}
\begin{itemize}
\item \emph{Robustness.}
Two sequences of graphs demonstrate that the assumption of robust connectivity is needed. In both these examples we use $\Uni[0,1]$-\previouslyweights\ on the edges. First,  take $G$ to be $K_{n-10}$ to which we attach a path $P_{10}$ of length 10. The minimum spanning tree on $G$ will obviously consist of a minimum spanning tree on $K_{n-10}$ plus all the edges of $P_{10}$. Thus the total \previouslyweight\ of the minimum spanning tree is in expectation~$5+o(1)$ more than that of the complete graph $K_n$, which $G$ is close to. Second, let $H$ and $H'$ consist of two cliques $A_1$ and $A_2$ of order $n/2$ each, which are either joined by a single complete vertex (in $H$) or a single edge (in $H'$). We see that the total \previouslyweight\ of the minimum spanning tree on $H$ is in expectation $0.5+o(1)$ less than in $H'$, even though these two graphs are close in the cut distance.

The latter example suggests that a possible relaxation of the notion of robustness is possible. Indeed, if instead of a single edge we connect $A_1$ and $A_2$ with, say, $n \log n$ edges then the total \previouslyweight\ of the minimum spanning tree on $H$ and $H'$ are very similar. The exact extent of this possible relaxation is not clear to us, as the problem of costly traversing across sparse cuts would reemerge in case of the graph being clustered into many moderately insulated parts (say more than $\log n$ in this case).
\item \emph{Decentness.}
Consider the sequence $(K_n)_n$ of cliques. We consider the distribution $D_n:=n^{-4}\cdot\Uni[0,n^{-4}]+(1-n^{-4})\cdot\Dirac(1)$ on each edge of $K_n$. We have that the derivative of each cumulative distribution function (on edges of $K_n$) is~1 at $0$. So, from this point of view, the setting is like in Theorem~\ref{thm:frieze}. In particular, if the (violated) condition~\ref{en:decent3} was not required then Theorem~\ref{thm:mainfull} would say that $(MST(G_n,D_n))_n$ converges in probability to $\zeta(3)$. This is however not the case in our setting. Indeed, for each $n$, all the edge \previouslyweights\ generated will be~$1$ with probability $1-\Theta(n^{-2})$. In that case, the minimum spanning tree has \previouslyweight\ $n-1\gg\zeta(3)$. By the Borel--Cantelli lemma this happens almost surely for all but finitely many $n$'s.
\end{itemize}

\subsection{Benjamini--Schramm convergence}
Our main result applies to sequences of dense graphs. Another well-developed graph limit theory is that of bounded degree graphs, in which the relevant topology is given by the Benjamini--Schramm convergence (also known as local or weak). It was noted in~\cite{MR2023650}, that if $(H_n)_n$ is a Benjamini--Schramm convergent sequence of connected graphs of growing orders and $D$ is a probability distribution with bounded support then $\frac{1}{v(H_n)}MST(H_n,D)$ converges in probability to a constant. An extension in which different distributions for different edges are used is also possible.\footnote{To that end, the notion of Benjamini--Schramm convergence needs to view edges as decorated by these distributions. Note also that the result only holds when all the supports of the distributions is uniformly bounded.} These results follow fairly easily by analysing Kruskal's algorithm. 

\subsection{Meaning of the constant $\kappa(W)$ via the branching process $\fX_W$}\label{ssec:meaningofconstant}
Recall that the giant component (or its nonexistence) in Erdős--Rényi random graphs $\RG(n,\frac{\lambda}{n})$ is typically studied via a Galton--Watson process whose offspring distributions are $\Poi(\lambda)$. In~\cite{MR2337396}, the problem of the giant component was studied in an inhomogeneous setting. To this end the following multitype branching process $\fX_W$ was introduced for each kernel $W\in L^\infty(\Omega^2)$. The first generation of $\fX_W$ consists of a single particle whose type has distribution $\mu$. Now, in any generation, a particle of type $x\in \Omega$ has offspring whose types are distributed as a Poisson process on $\Omega$ with intensity $W(x,y)\differential\mu(y)$. That is, the number of offspring whose types are in a given set $Z\subset \Omega$ is $\Poi(\int_{y\in Z}W(x,y)\differential\mu(y))$, and for disjoint sets $Z_1$ and $Z_2$ the respective numbers of offspring of these types are independent. The process $\fX_W$ played a key role also in~\cite{MR2599196,MR2880662,MR2816939}. Note that when a realization of $\fX_W$ has finite progeny, it can be represented as a \emph{plane tree}, that is, a finite or an infinite tree which is rooted and where children of each vertex are linearly ordered.

We can now define the key quantities for Theorem~\ref{thm:mainfull}. For a kernel $W$, define
\begin{align}
\label{eq:defUpsilon}
\Upsilon(W)&:=\sum_{k=1}^\infty \frac{\Pr[|\fX_{W}|=k]}k \; \mbox{and}\\
\label{eq:defkappa}
\kappa(W)&:=\int_{c=0}^{+\infty}\Upsilon(c\cdot W)\differential c\;.
\end{align}
There are kernels $W$ for which $\kappa(W)=+\infty$. The simplest case is the constant-0 kernel. Indeed, $\fX_0$ almost surely consists only of the root. Thus, $\Upsilon(0)=\frac{1}{1}$, and so $\kappa(0)=\int_{c=0}^{+\infty}\Upsilon(c\cdot 0)=+\infty$. This argument generalizes easily to kernels $W$  in which $\deg_W(x)=0$ for a positive measure of $x\in\Omega$, and its refined version generalizes even to some kernels $W$ in which $\deg_W(x)>0$ for every $x\in \Omega$ but $\mindeg(W)=0$.
Theorem~\ref{thm:mainfull} claims that $\kappa(W)<\infty$ for kernels that arise as cut distance limits of graphs $(G_n)_n$ from Setup~\ref{mainsetup}. The two facts below, whose proofs we defer to Section~\ref{sec:proofsKappa}, justify this.
\begin{fact}\label{fact:robustmindeg}
	Suppose that we have $\rho>0$ and weighted graphs $(G_n)_n$ which are all $\rho$-robust. If a kernel $W$ is a cut distance limit of $(G_n)_n$ then $\mindeg(W)\ge \frac{\rho}2$.
\end{fact}
\begin{fact}\label{fact:kappafinite}
	Suppose that $W$ is a kernel with a positive minimum degree. Then $\kappa(W)<\infty$.
\end{fact}

We now give a definition of a fractional multiples which allows us to relate the parameter $\kappa$ for some pairs of kernels. This notion is new and our choice of the term comes from the fact that it extends the concept of `fractional isomorphism' of graphons worked out in~\cite{MR4482093}, which corresponds to 1-fractional multiples. After the first version of this paper was posted on arXiv, the same concept appeared in~\cite{HlaHngLim:GraphonBranching} under the name `projective fractional isomorphism'. 
\begin{definition}\label{def:fractionalmultiple}
For a constant $p>0$, we say that a kernel $U$ is a \emph{$p$-fractional multiple} of a kernel $W$, if for each tree $T$, we have for the homomorphism densities (recall~\eqref{eq:homdens}) that
\begin{equation}\label{eq:TreeDens}
t(T,U)=p^{e(T)}t(T,W)    \;.
\end{equation}
\end{definition}
The most prominent example of fractional multiples are regular kernels. If $U$ is an arbitrary $d$-regular kernel and $W$ is an arbitrary $q$-regular kernel then a quick calculation shows that $t(T,U)=d^{e(T)}$ and $t(T,W)=q^{e(T)}$ for each tree~$T$. We conclude that $U$ is a $\frac{d}{q}$-fractional multiple of $W$. 

The primary contribution of~\cite{MR4482093} lies in presenting several characterizations of fractional isomorphism, all of which naturally extend to characterizations of $p$-fractional multiples. The significance of this concept is highlighted by the following fact.
\begin{fact}\label{fact:AFM}
Suppose that $U$ and $W$ are two kernels such that $U$ is $f$-fractional multiple of $W$ for some $f>0$. Then $\kappa(U)=\frac{1}{f}\cdot\kappa(W)$.
\end{fact}
We deduce this fact from one of the main results of~\cite{HlaHngLim:GraphonBranching}. Indeed, \cite{HlaHngLim:GraphonBranching} tells that if $A$ and $B$ are two kernels so that $A$ is a 1-fractional multiple of $B$, then the branching processes $\fX_A$ and $\fX_B$ have the same distribution. In particular, $\Upsilon(A)=\Upsilon(B)$. The substitution formula for integration with $\tilde{c}=f c$ (and hence $\differential\tilde{c}=f \differential{c}$) in~\eqref{eq:defkappa} then readily implies that $\kappa(U)=\frac{1}{f}\cdot\kappa(W)$.

Last, let us look at the most prominent kernel, the constant-1. We have $\kappa(\mathbbm{1})=\zeta(3)$. Indeed, complete graphs (of growing orders) equipped with $\Uni[0,1]$-\previouslyweights\ converge to the kernel $\mathbbm{1}$ in the cut distance. By Theorem~\ref{thm:mainfull}, the total \previouslyweights\ of the random minimum spanning trees converge to a constant $\kappa(\mathbbm{1})$ in probability and by Theorem~\ref{thm:frieze}, they converge to~$\zeta(3)$ in probability, so $\kappa(\mathbbm{1})=\zeta(3)$. This argument is obviously not self-contained. In Section~\ref{ssec:ComputingZeta3} we give a self-contained proof inspired by Frieze's calculations.

\subsubsection{Theorem~\ref{thm:mainfull} versus Theorem~\ref{thm:FriezeMcDiarmid}}\label{ssec:MAINvsFrMcD}
Let us show that Theorem~\ref{thm:FriezeMcDiarmid} is implied by Theorem~\ref{thm:mainfull}. Suppose that $\mathfrak{H}=(H,(D_e)_{e\in E(H)})$ is a $\Delta$-regular template in which $H$ is a connected graph. Then it is pedestrian to show that the sequence $(G_n,\cD_n)$ of $n$-blow-ups of $\mathfrak{H}$ is decent (there are only finitely many distributions involved). Further, the weighted graphs $G_n$ corresponding to $(H_n,\cD_n)$ are all $\rho$-robust where $\rho$ can be taken as the minimum of the weights of~$H$ (recall that in the definition of regular templates, all the derivatives of the cumulative distribution functions are positive) divided by~$v(H)$. Last, the sequence $(G_n)_n$ converges to a kernel $W$ which is a block kernel consisting of square of $v(H)$-many blocks in $\Omega$ of measure $\frac{1}{v(H)}$ each, where the value of $W$ on a product of blocks corresponding to vertices $u$ and $v$ of $H$ is the value of derivative of the cumulative distribution function of $D_{uv}$ at~0 if $uv\in E(H)$ and~0 otherwise. Hence, $\Delta$-regularity of $\mathfrak{H}$ translates as $\Delta$-regularity of $W$. It is easy to check that then $W$ is a $\Delta$-fractional multiple of the constant-1 kernel $\mathbbm{1}$, and hence Theorem~\ref{thm:FriezeMcDiarmid} follows as a consequence.

On the other hand, the main new feature (apart from the fact that it can handle infinitesimal setting thanks to the formalism of graph limits) compared to Theorem~\ref{thm:FriezeMcDiarmid} is that the degree regularity is not needed.

\subsubsection{Some examples} 
\begin{figure}
	\includegraphics[scale=0.75]{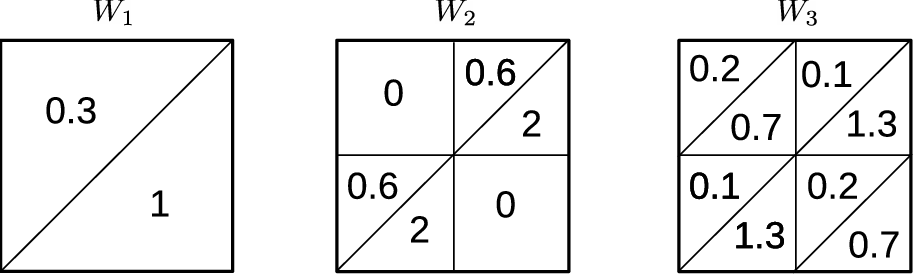}
	\caption{Three fractionally isomorphic graphons $W_1$, $W_2$, and $W_3$. (In the figures, we use the same orientation of the plane as for matrices, that is, the main diagonal starts in the top left corner.)}
	\label{fig:three}
\end{figure}

In view of Section~\ref{ssec:MAINvsFrMcD} we focus on nonregular examples. In Figure~\ref{fig:three} we show three graphons $W_1$, $W_2$, and $W_3$ which are fractionally isomorphic, that is, they are 1-fractional multiples of each other. This fact can be verified using~\eqref{eq:TreeDens}, but equivalent characterizations given in~\cite{MR4482093} are even more convenient to this end. Hence, we have $\kappa(W_1)=\kappa(W_2)=\kappa(W_3)=:\kappa^*$.
Theorem~\ref{thm:mainfull} tells us that for any large graph (satisfying the technical assumptions of Theorem~\ref{thm:mainfull}) equipped with distributions on the edges whose corresponding weighted graph is close to any of these graphons, the \previouslyweight\ of the random minimum spanning tree is concentrated close to $\kappa^*$.

Let us give some examples for these three graphons.
\begin{itemize}
	\item Consider a complete graph on vertex set $[n]$. Each edge $ij$ is equipped with a distribution $\Uni[0,1]$ if $i+j>n$ and with a distribution $\Uni[0,1/0.3]$ if $i+j\le n$. The corresponding weighted graph is close to $W_1$.
	\item Consider a graph on vertex set $[n]$. All the edges are equipped with the distribution $\Uni[0,1]$. Pairs $ij$ are inserted as edges with probability $1$ if $i+j>n$ and with probability $0.3$ if $i+j\le n$. The corresponding weighted graph typically is close to $W_1$, too.
	\item Consider a complete bipartite graph on vertex set $\{1,2,\ldots,n\}$ where the two parts are formed by even and odd vertices, respectively. On edges $ij$ ($i$ is even, $j$ is odd) put a distribution $\Uni[0,1/2]$ if $i+j>n$ and $\Uni[0,1/0.6]$ if $i+j\le n$. The corresponding weighted graph is close $W_2$. A similar construction could be given for $W_3$.
\end{itemize}
Section~2 of~\cite{MR4482093} gives a general (and exhaustive) approach for constructing families of fractionally isomorphic graphons. This gives many more further examples.

\section{Proof of Theorem~\ref{thm:mainfull}}
Our proof of Theorem~\ref{thm:mainfull} relies on three auxiliary lemmas which all deal with connectivity properties of certain random graphs. Proofs of these lemmas are left to subsequent sections. Let us introduce notation necessary to state these lemmas. Suppose that $G$ is a graph of order $n$ equipped with a weight function $w:E(G)\to [0,\infty)$. Then the \emph{percolation} of $G$, denoted $\mathsf{percolation}(G)$, is a random (unweighted) subgraph of $G$ in which we keep each edge $e\in E(G)$ with probability $\min(1,\frac{w(e)}{n})$ and all the choices are independent. Suppose that $H$ is an $n$-vertex graph. Write $cc(H)$ for the number of connected components of $H$. Then we call the quantity $dr(H):=\frac{cc(H)}{n}$ the \emph{disconnectedness ratio} of $H$.

The first auxiliary lemma says that if we have two weighted graphs that are close in the cut distance then the disconnectedness ratios of their percolations are with high probability close. Furthermore, the disconnectedness ratio relates to the branching process $\fX_\Gamma$ (of any kernel $\Gamma$ which is close to those graphs in the cut distance).
\begin{lemma}\label{lem:FirstLemma}
Suppose that $(G_n)_n$ is a sequence of weighted graphs of growing orders with uniformly bounded weights that converge to a kernel $\Gamma$ in the cut distance. Then the disconnectedness ratios of the percolations of $G_n$ converge in probability to the number $\Upsilon(\Gamma)$ defined in~\eqref{eq:defUpsilon}.
\end{lemma}
The next lemma asserts that a random subgraph of a robust $n$-vertex graph in which edges are kept with probability $\Theta(\frac{\log n}{n})$ is connected with high probability.
\begin{lemma}\label{lem:SecondLemma}
For every $\alpha>0$, the following holds. Suppose that $H$ is an $\alpha$-robust $n$-vertex graph. Let $F\subset H$ be a random subgraph of $H$ in which each edge is included independently with probability $\frac{10\log n}{\alpha n}$. Then with probability at least $1-2 n^{-3}$, the graph $F$ is connected.
\end{lemma}
Note that Lemma~\ref{lem:SecondLemma} is related to~\cite{DevFra:Connectivity}, which studies connectivity of inhomogeneous random graph models where edge probabilities scale as $\frac{\log n}{n}$. We cannot use~\cite{DevFra:Connectivity} directly, but our proof is short anyway (mostly because the involved constant $\frac{10}{\alpha}$ is not optimal). 

Our last auxiliary lemma asserts that if we add to an $n$-vertex graph $F$ random edges with success probability $\Theta(\frac{1}{n})$ from a robust template $H$, the number of connected components in the resulting graph drops substantially with high probability. 
\begin{lemma}\label{lem:CompCountDrop}
Suppose that $\alpha>0$. Suppose that $F$ and $H$ are two graphs on the same vertex set $V$ with $|V|=n$ satisfying $n>\frac{500}{\alpha}$. Suppose that $H$ is $\alpha$-robust. Let $F'$ be obtained from $F$ by adding independently edges of $E(H)\setminus E(F)$, so that an individual edge from $E(H)\setminus E(F)$ is added with probability $\frac{12}{\alpha n}$. Then with probability at least $1-n^{-4}$, we have $cc(F')\le 0.6 cc(F)+n^{4/5}$.
\end{lemma}
Lemma~\ref{lem:FirstLemma} is proven in Section~\ref{sec:Proof1}, Lemma~\ref{lem:SecondLemma} is proven in Section~\ref{sec:Proof2} and Lemma~\ref{lem:CompCountDrop} is proven in Section~\ref{sec:Proof3}.

\medskip

As with similar results in the area, our proof of Theorem~\ref{thm:mainfull} relies on the analysis of Kruskal's algorithm. Recall that given a connected graph $F$ with a \previouslyweight\ function $w$, starting from the edgeless graph on $V(F)$, in each step of Kruskal's algorithm,
we add an edge of the smallest \previouslyweight\ that will not introduce a cycle to the current forest (if there are several edges of the smallest \previouslyweight\ the algorithm can make an arbitrary choice). A crucial property of Kruskal's algorithm is that it produces a minimum spanning tree. 
For $x\ge 0$, we denote by $F^{\le x}\subset F$, the unweighted spanning subgraph of $F$ consisting of edges of \previouslyweight\ at most $x$. Suppose that $T$ is a minimum spanning tree of $F$. Then for each $e\in E(T)$ and for each
$x<w(e)$, the two endvertices of $e$ lie in different components of the graph $F^{\le x}$, while for each $x\ge w(e)$, they lie in the same component of the graph $F^{\le x}$. This gives the following identity,
\begin{equation}\label{eq:MSTNumberComps}
MST(F)=\int_{x=0}^\infty\left(cc\big(F^{\le x}\big)-1\right)\differential x=\int_{x=0}^\infty\big(v(F)\cdot dr\big(F^{\le x}\big)\underbrace{-1}_{\textsf{(T)}}\big)\differential x\;.
\end{equation}
This identity (even if only implicit) was crucial in previous work on the random minimum spanning tree problem.

We can now proceed to the main part of the proof of Theorem~\ref{thm:mainfull}. 
Suppose that $\gamma>0$ is arbitrary.
We set some key constants for the proof. We take $\Delta:=\gamma/4$. Further, we take $\eps>0$, such that 
\begin{equation}\label{eq:epsilons}
	\kappa\big((1-\eps)W\big)\le \kappa(W)+\frac\gamma2 \quad\mbox{and}\quad
	\kappa(W)\le \kappa\big((1+\eps)W\big)+\frac\gamma2 \;.
\end{equation}
Lemma~\ref{lem:MultiContinuityofKappa} tells us that this is possible.

Last, we fix a constant $K\in \NN$ so large that we have the following two conditions. Firstly, we require that $\int_{x=\Delta K}^{\infty}(\rho x)^{-2}<\gamma/40$. In particular, Lemma~\ref{lem:kappafinitedeg}\ref{en:kappaUpsilon} tells us that 
\begin{equation}\label{eq:Ktail}
	\int_{c=\Delta K}^{\infty}\Upsilon(cW)\differential c<\frac\gamma8\;.
\end{equation}
Secondly, we require that for every kernel $Z$ with minimum degree at least $\frac{\rho^2\Delta}{6}$, we have 
\begin{equation}\label{eq:KKK}
\Upsilon(KZ)\le \frac{\gamma\cdot \rho^2}{20000}\;.
\end{equation}
Lemma~\ref{lem:kappafinitedeg}\ref{en:kappaUpsilon} tells us that this is possible.

For simplicity and without loss of generality, assume that each $H_n$ has $n$ vertices.
We shall think of the graphs $H_n$ as equipped with their random \previouslyweights\ $\cD_n$. To refer to this graph with random \previouslyweights\ of the edges, we write $\G_n$. In particular, $MST(\G_n)$ is a random variable and $\G_n^{\le x}$ is a random graph (for each $x\ge 0$). We approximate the integral in~\eqref{eq:MSTNumberComps} by Riemann sums with step size $\frac{\Delta}{n}$. We use monotonicity of $dr(F^{\le x})$ in the variable $x$. The upper bound and the lower bound below hold for an arbitrary $K\in\NN$, and we shall fix it as a sufficiently large constant later. For the upper bound, we split the integration~\eqref{eq:MSTNumberComps} at $x=\frac{K\Delta}{n}$, and also omit the term $\textsf{(T)}$ appearing in~\eqref{eq:MSTNumberComps}.
\begin{align}
\label{eq:KeyLowerBound}
MST(\G_n)&\le \Delta\sum_{i=0}^{K-1} dr\left(\G_n^{\le \frac{i\Delta}n}\right)+\int_{x=\frac{K\Delta}{n}}^\infty\left(cc\big(\G_n^{\le x}\big)-1\right)\differential x\;.
\end{align}
For the lower bound, we remove the tail of the integral~\eqref{eq:MSTNumberComps} corresponding to $x\ge \frac{K\Delta}{n}$. Note that the contribution of the term~$\textsf{(T)}$ in~\eqref{eq:MSTNumberComps} to the remaining (retained) portion is $-\frac{K\Delta}{n}$.
\begin{align}
\label{eq:KeyUpperBound}
MST(\G_n)&\ge \Delta\sum_{i=0}^{K-1} dr\left(\G_n^{\le \frac{(i+1)\Delta}n}\right) -\frac{K\Delta}{n}\;.
\end{align}
We first analyze the sums involving the random graphs $\G_n^{\le \frac{j\Delta}n}$ in~\eqref{eq:KeyLowerBound} and~\eqref{eq:KeyUpperBound}. A particular edge $e\in E(G_n)$ is contained in $\G_n^{\le \frac{j\Delta}n}$ with probability $D_{n,e}([0,\frac{j\Delta}n])$. We use~\ref{en:decent3} of Definition~\ref{def:decent} to see that $D_{n,e}([0,\frac{j\Delta}n])=(1\pm \eps)F_{n,e}'(0)\frac{j\Delta}n$ (when $j $ and $\Delta$ are bounded by constants and $n$ is sufficiently large). For each $j=0,\ldots,K$, we define weighted graphs $G_{n}^{j,-}:=(1-\eps)(j\Delta)G_n$ and $G_{n}^{j,+}:=(1+\eps)(j\Delta)G_n$, by which we mean that we multiply by constants $(1-\eps)(j\Delta)$ or $(1+\eps)(j\Delta)$, the weights of the graph $G_n$, which was defined to be the weighted graph corresponding to $H_n$ and $\cD_n$. 

We know that the sequence $(G_n)_n$ converges to $W$. Hence, for each $j=0,\ldots,K$, the sequence $(G_{n}^{j,-})_n$ converges to $(1-\eps)(j\Delta)W$ and the sequence $(G_{n}^{j,+})_n$ converges to $(1+\eps)(j\Delta)W$. In both these sequences, the weights are uniformly bounded. So, we are in the setting of Lemma~\ref{lem:FirstLemma}.  Hence for $j=0,\ldots,K$, we have
\begin{align}
	\label{eq:MLower}
dr\left(\G_n^{\le \frac{j\Delta}n}\right) &\leBy{stoch domination} dr\left(\mathsf{percolation}(G_{n}^{j,-})\right) \stackrel{n\to\infty}{\longrightarrow}
\Upsilon\big((1-\eps)(j\Delta)W\big)
\mbox{\;, and}\\
	\label{eq:MUpper}
dr\left(\G_n^{\le \frac{j\Delta}n}\right) &\geBy{stoch domination} dr\left(\mathsf{percolation}(G_{n}^{j,+})\right) \stackrel{n\to\infty}{\longrightarrow}
\Upsilon\big((1+\eps)(j\Delta)W\big)\;.
\end{align}

We can now analyze the (random) term 
\begin{equation}\label{eq:wli}
\int_{x=\frac{K\Delta}{n}}^\infty\left(cc\big(\G_n^{\le x}\big)-1\right)\differential x
\end{equation}
in~\eqref{eq:KeyLowerBound}. Let $R_n$ be the unweighted version of $G_n$ in which only edges of weight at least $\frac\rho2$ are kept. Also, let $M$ be the supremum of all the edge weights (over all the graphs $\{G_n\}_n$; recall~\eqref{eq:sup}).
\begin{claim}
Each graph $R_n$ is $\rho/(2M)$-robust.
\end{claim}
\begin{proof}
Let $E_n^-\subset E(G_n)$ be the edges of $G_n$ of weight less than $\frac\rho2$. In the calculation below, we write $e_{G_n\setminus E_n^-}(\cdot,\cdot)$ for the total weight of the edges between two sets in the graph $G_n\setminus E_n^-$. Let $U\subset V(G_n)$ be arbitrary. Using that $G_n$ is $\rho$-robust, we have
\[
e_{G_n\setminus E_n^-}(U,V(G_n)\setminus U)\ge e_{G_n}(U,V(G_n)\setminus U)-\tfrac\rho2 |U|(n-|U|)\ge \tfrac\rho2 |U|(n-|U|)\;.
\]
In particular, since the weight of each edge of $G_n$ is at most $M$, we get for the unweighted graph $R_n$ that $e_{R_n}(U,V(G_n)\setminus U)\ge \tfrac\rho{2M} |U|(n-|U|)$, as was needed.
\end{proof}
For $p\in[0,1]$, let $\mathbf{R}_n(p)$ be a random subgraph of $R_n$ in which each edge is kept with probability $p$. We claim that for $n$ sufficiently large, for each $p\le \frac{20\log n}{\rho n}$,
\begin{equation}\label{eq:graphdomin}
\mbox{the graph $\mathbf{R}_n(p)$ is stochastically dominated by the graph $\G_n^{\le \frac{3p}\rho}$.}
\end{equation}To see this, we consider a particular edge $e\in E(H_n)$ and distinguish whether its weight is less than or at least $\frac{\rho}{2}$. In the former case, $e$ appears with probability~0 in $\mathbf{R}_n(p)$, so domination on the edge~$e$ is trivial. In the latter case, $e$ appears with probability~$p$ in $\mathbf{R}_n(p)$ and it appears with probability 
$$D_{n,e}\left(\left[0,\tfrac{3p}\rho\right]\right)=(1+o(1))\cdot\frac{3p}\rho\cdot F_{n,e}'(0)\ge (1+o(1))\cdot\frac{3p}\rho\cdot \frac{\rho}{2}\ge p$$
in $\G_n^{\le \frac{3p}\rho}$, as was needed for~\eqref{eq:graphdomin}. In particular, $cc(\G_n^{\le 3p/\rho})$ is stochastically dominated by $cc(\mathbf{R}_n(p))$.

We split the integration in~\eqref{eq:wli} into two domains, namely the interval $\left[\frac{K\Delta}{n},\frac{60\log n}{\rho^2 n}\right]$ and the interval $\left[\frac{60\log n}{\rho^2 n},\infty\right)$.

First, we deal with integration on the interval $\left[\frac{K\Delta}{n},\frac{60\log n}{\rho^2 n}\right]$ in~\eqref{eq:wli}. To this end, we observe that for any $p,q\in[0,1]$ with $p+q\le 1$, the random graph $\mathbf{R}_n(p+q)$ stochastically dominates the edge-union of the random graphs $\mathbf{R}_n(p)$ and $\mathbf{R}_n(q)$. We use this repeatedly with $$p_i:=\frac{\rho}{3}\cdot\frac{K\Delta}{n}+i\cdot \frac{12}{\frac{\rho}{2M}n}$$ (where $i=0,\ldots,\lfloor n^{0.01}\rfloor$) and a universal $q:=\frac{12}{\frac{\rho}{2M}n}$. For each such $i$, Lemma~\ref{lem:CompCountDrop} tells us that, up to an error probability of at most $n^{-4}$, we have $cc(\mathbf{R}_n(p_i)\cup \mathbf{R}_n(q))\le 0.6cc(\mathbf{R}_n(p_i))+ n^{4/5}$. In particular, by the above stochastic domination, we have $cc(\mathbf{R}_n(p_{i+1}))\le 0.6cc(\mathbf{R}_n(p_i))+ n^{4/5}$, except for an event of probability at most $n^{-4}$. Using the union bound, we have that with probability at least $1-n^{-3.99}$, it holds that for each $\ell=1,\ldots,\lfloor n^{0.01}\rfloor$, we have
\begin{equation}\label{eq:Am}
cc(\mathbf{R}_n(p_{\ell}))	\le 0.6^\ell cc(\mathbf{R}_n(p_{0})) + \sum_{i=0}^{\ell-1} 0.6^i n^{4/5}\le 0.6^\ell cc(\mathbf{R}_n(p_{0})) + 2.5 n^{4/5}\;.
\end{equation}
Fix $L:=\left(\frac{60\log n}{\rho^2 n}-\frac{K\Delta}{n}\right)/(\rho/3)$. Therefore, with probability at least $1-n^{-3.99}$, Riemann summation with steps of length $\frac{3q}{\rho}$ gives
\begin{align}
\begin{split}
\label{eq:integraltail1}
	\int_{x=\frac{K\Delta}{n}}^{\frac{60\log n}{\rho^2 n}} \left(cc\big(\G_n^{\le x}\big)-1\right)\differential x
	&\leByRef{eq:graphdomin}  
	\frac{3q}{\rho}\cdot\sum_{\ell=0}^{L}\left(cc(\mathbf{R}_n(p_{\ell}))-1\right)\\
\JUSTIFY{by~\eqref{eq:Am}}&\le \frac{3q}{\rho}\cdot\sum_{\ell=0}^{\infty}0.6^\ell cc(\mathbf{R}_n(p_{0})) \; + \; \frac{3q}{\rho}\cdot L\cdot 2.5 n^{4/5}\\
\JUSTIFY{$q\ll n^{-0.99}$, $L\ll n^{0.01}$}&\le \frac{3q}{\rho}\cdot\sum_{\ell=0}^{\infty}0.6^\ell cc(\mathbf{R}_n(p_{0})) \; + \; n^{-0.18}\\
&=
\frac{7.5q}{\rho}\cdot cc(\mathbf{R}_n(p_{0})) +  n^{-0.18}\;.
\end{split}
\end{align}
We apply~\eqref{eq:KKK} to the kernel $Z$ which is the cut distance limit of the sequence $(\frac{\rho}{3}\cdot\Delta \cdot R_n)_n$, which indeed satisfies the above minimum-degree condition.\footnote{Strictly speaking, the sequence $(R_n)_n$ need not be convergent in the cut distance. So, should the calculations below fail for a subsequence $(R_{n_i})_{n_i}$ we could use compactness of the cut distance topology (Theorem~\ref{thm:compactness}), find 
a convergent subsequence $(R_{m_i})_{m_i}$ and a limit kernel $Z$. Then the calculations below show that the sequence $(R_{n_i})_{n_i}$ was not a sequence of counterexamples.}
 Lemma~\ref{lem:FirstLemma} tells us that $dr(\mathbf{R}_n(p_{0}))$ converges in probability to $\Upsilon(KZ)\le \frac{\gamma\cdot \rho^2}{20000}$. In particular, $\frac{7.5q}{\rho}\cdot cc(\mathbf{R}_n(p_{0}))$ converges in probability to a constant which is at most $\frac{\gamma}{100}$. Substituting to~\eqref{eq:integraltail1}, we get that with probability at least $1-\frac{\gamma}{10}$, for large enough $n$ we have,
\begin{equation}\label{eq:integraltailA}
	\int_{x=\frac{K\Delta}{n}}^{\frac{60\log n}{\rho^2 n}} \left(cc\big(\G_n^{\le x}\big)-1\right)\differential x\le \frac{\gamma}{99}\;.
\end{equation}

We now deal with integration on the interval $\left[\frac{60\log n}{\rho^2 n},\infty\right)$ in~\eqref{eq:wli}. Lemma~\ref{lem:SecondLemma} asserts that asymptotically almost surely, $\mathbf{R}_n(\frac{20\log n}{\rho n})$ is connected. Put together with~\eqref{eq:graphdomin}, asymptotically almost surely $\G_n^{\le 60\log n/(\rho^2 n)}$ is connected. In particular, asymptotically almost surely, 
\begin{equation}\label{eq:integraltail2}
\int_{x=\frac{60\log n}{\rho^2 n}}^\infty \left(cc\big(\G_n^{\le x}\big)-1\right)\differential x=0\;.
\end{equation}

We now plug the above bounds into~\eqref{eq:KeyLowerBound} and~\eqref{eq:KeyUpperBound}. 

Let us start with~\eqref{eq:KeyLowerBound}. We use~\eqref{eq:MLower} for the main term and~\eqref{eq:integraltailA} and~\eqref{eq:integraltail2} for the tail integral. We get that for $n$ sufficiently large with probability at least $1-\frac{\gamma}9$ we have,
\begin{align*}
MST(\G_n)&\le \Delta\sum_{i=0}^{K-1} 
\Upsilon\big((1-\eps)(i\Delta)W\big)
\;
+\;
\frac{\gamma}{99}
\;
+
\;
0\;.
\end{align*}
Since $\Upsilon(\cdot)\le 1$ and thanks to the monotonicity property (Lemma~\ref{lem:UpsilonMonot}), we can view the sum as a Riemann integration. Thus, with probability at least $1-\frac{\gamma}9$, 
\begin{align}\label{eq:JB1}
	MST(\G_n)&\le \Delta+	\frac{\gamma}{99}+\int_0^{K\Delta}
	\Upsilon\big((1-\eps)c W\big)\differential c
	\le 
	\Delta+	\frac{\gamma}{99}+\kappa\big((1-\eps) W\big)
	\leByRef{eq:epsilons} 
	\gamma+\kappa(W)\;.
\end{align}

The lower bound follows similarly. That is, we get from~\eqref{eq:KeyUpperBound} and~\eqref{eq:MUpper} that with probability at least $1-\frac{\gamma}9$, 
\begin{align*}
	\frac{K\Delta}{n}+MST(\G_n)&\ge \Delta\sum_{i=1}^{K} 
	\Upsilon\big((1+\eps)(i\Delta)W\big)\ge \int_0^{K\Delta}
	\Upsilon\big((1+\eps)c W\big)\differential c \\
	&=\kappa\big((1+\eps)W\big)-\int_{K\Delta}^\infty
	\Upsilon\big((1+\eps)c W\big)\differential c\\
	\JUSTIFY{substitution $c^*=(1+\eps)c$}&=\kappa\big((1+\eps)W\big)-\frac{1}{1+\eps}\cdot\int_{(1+\eps)K\Delta}^\infty
	\Upsilon\big(c^* W\big)\differential c^*
	  \\
	  \JUSTIFY{by \eqref{eq:epsilons} and \eqref{eq:Ktail}}&\ge  \kappa(W)-\frac\gamma2-\frac{1}{1+\eps}\cdot\frac\gamma8\;.
\end{align*}
Since $K$ and $\Delta$ were constants independent of $n$, the term $\frac{K\Delta}{n}$ vanishes as $n\to\infty$. Hence, with probability at least $1-\frac{\gamma}9$, 
\begin{equation}\label{eq:JB2}
    MST(\G_n)\ge \kappa(W)-\gamma\;.
\end{equation}
Since $\gamma>0$ was arbitrary, the theorem follows from~\eqref{eq:JB1} and~\eqref{eq:JB2}.

\subsection{Proof of Lemma~\ref{lem:FirstLemma}}\label{sec:Proof1}
Lemma~\ref{lem:FirstLemma} follows quite easily from the main lemma of~\cite{MR2599196}. That lemma talks about the number of components of a specific size in percolations of $G_n$. To this end, in an $n$-vertex graph~$H$, we write $\ell_k(H)$ for the number of connected components of~$H$ of order exactly~$k$. The number $ir_k(H)=\frac{\ell_k(H)}n$ is the \emph{$k$-island ratio} of~$H$.
\begin{lemma}[Lemma~4 in~\cite{MR2599196}]\label{lem:BBCR}
Suppose that $(G_n)_n$ is a sequence of weighted graphs with uniformly bounded weights that converge to a kernel $\Gamma$ in the cut distance. Then for each $k\in \NN$, the $k$-island ratios of the percolations of $G_n$ converge in probability to $\frac{\Pr[|\fX_{\Gamma}|=k]}k$.
\end{lemma}

Suppose that $F$ is an arbitrary graph. Then we have
\begin{equation*}
dr(F)=\sum_{k=1}^\infty  ir_k(F)\;.
\end{equation*}
We want to obtain bounds on $dr(F)$ by discarding large components, say of order more than $L\in\NN$. The number of such components is between~0 and $\frac{v(F)}{L+1}$. Hence we arrive at
\begin{equation}\label{eq:reckless}
\frac{1}{L+1}+\sum_{k=1}^L  ir_k(F) \ge dr(F)\ge \sum_{k=1}^L  ir_k(F)\;.
\end{equation}

With these preparations, we can prove Lemma~\ref{lem:FirstLemma} readily. Let $\eps>0$ be arbitrary. Take $L:=\lceil\frac2\eps\rceil$. Lemma~\ref{lem:BBCR} tells us that asymptotically almost surely, for each $k\in [L]$, the $k$-island ratios of the percolations of $G_n$ are in the interval $\left[\frac{\Pr[|\fX_{\Gamma}|=k]}k-\frac{\eps}{2L},\frac{\Pr[|\fX_{\Gamma}|=k]}k+\frac{\eps}{2L}\right]$. Using~\eqref{eq:reckless}, we obtain that the disconnectedness ratio of the percolations of $G_n$ is asymptotically almost surely in the interval whose minimum is $$\sum_{k=1}^L\left(\frac{\Pr[|\fX_{\Gamma}|=k]}k-\frac{\eps}{2L}\right)\ge \sum_{k=1}^L\frac{\Pr[|\fX_{\Gamma}|=k]}k-\frac\eps2\ge \sum_{k=1}^\infty\frac{\Pr[|\fX_{\Gamma}|=k]}k-\eps$$ 
and whose maximum is $$\frac{1}{L+1}+\sum_{k=1}^L\left(\frac{\Pr[|\fX_{\Gamma}|=k]}k+\frac{\eps}{2L}\right)\le \sum_{k=1}^\infty\frac{\Pr[|\fX_{\Gamma}|=k]}k+\eps\;.$$
As $\eps$ was arbitrary, this finishes the proof.

\subsection{Proof of Lemma~\ref{lem:SecondLemma}}\label{sec:Proof2}
We have,
\begin{align*}
\Pr\left[\text{$F$ is disconnected}\right] &\leq \sum_{U \subseteq V(H), |U| \leq \frac{n}2, U\neq \emptyset} \left( 1-\frac{10 \log n}{\alpha n} \right)^{\alpha |U| (n-|U|)} \\
&=  \sum_{k=1}^{\lfloor \frac{n}2\rfloor} 
\binom{n}{k} \cdot\left( 1-\frac{10 \log n}{\alpha n} \right)^{\alpha k (n-k)}\leq \sum_{k=1}^{\lfloor \frac{n}2\rfloor} 
n^{k} \cdot \left( 1-\frac{10 \log n}{\alpha n} \right)^{\frac{\alpha k n}{2}}\\
&\le \sum_{k=1}^{\lfloor \frac{n}2\rfloor}  \exp(k \log n -5 k \log n)\leq \frac{n}{2} \cdot\exp(-4 \log n) \leq \frac{2}{ n^{3}}\;.
\end{align*}
\subsection{Proof of Lemma~\ref{lem:CompCountDrop}}\label{sec:Proof3}
Call a component of $F$ \emph{untouched} if it is also a component of $F'$. Let $r$ be the number of untouched components.

We need a counterpart of robustness for oriented graphs. For $\sigma>0$, we say that an $n$-vertex oriented graph $\overrightarrow{G}$ is \emph{$\sigma$-tough}, if for every $U\subset V(\overrightarrow{G})$ , the number of oriented edges of $\overrightarrow{G}$ from $U$ to $V(\overrightarrow{G})\setminus U$ is at least $\sigma |U|(n-|U|)$. 
\begin{claim}\label{cl:orbexist}
There exist an $\frac{\alpha}{3}$-tough orientation 
$\overrightarrow{H}$ of $H$.
\end{claim}
\begin{proof}
We randomly orient each edge of $H$ with probability $1/2$ in each direction. 
Fix  $U \subseteq{V(H)}$, $|U| \leq \frac{n }{2}$. Let $X$ denote the random variable which is number of edges oriented from $U$ to $V(H)\setminus U$.  Then $X$ can be written as a sum of $|U|$ independent random variables, where each random variable represents the number of oriented edges from a vertex in $U$.
Let $E_{U}$ be the event that there are less than $ \frac{\alpha |U| (n-|U|)}{3}$ edges oriented from $U$ to $V(H)\setminus U$.
Since $H$ is $\alpha$-robust, $e_{H}(U,V(H)\setminus U) \geq \alpha |U| (n-|U|)$. Hence $\Exp[X] \geq \frac{\alpha |U| (n-|U|)}{2} $. By Chernoff's bound,
\[
\Pr\left[E_{U}\right] \leq \Pr \left[X \leq \left(1-\frac{1}{3} \right) \Exp[X] \right] \leq \exp(-\Exp[X]/18) \leq \exp(-\alpha |U| (n-|U|)/36 ).
\]
We shall show that  $\sum_{U  \subseteq V(H) , \ |U| \leq n /2, U\neq \emptyset} \Pr\left[E_{U}\right]< 1 $ and the claim will follow by the probabilistic method.
Indeed, we have
\begin{align*}
\sum_{U  \subseteq V(H), \ |U| \leq n /2, U\neq \emptyset } \Pr\left[E_{U}\right] &\leq \sum_{k=1}^{\lceil\frac{n}2\rceil} \binom{n}{k} \exp(-\alpha k (n-k)/36 )\\ 
&\leq \sum_{k=1}^{\lceil\frac{n}2\rceil} n^{k} \exp(-\alpha k n/73 )= \sum_{k=1}^{\lceil\frac{n}2\rceil}  \exp\big(k (\log n-\alpha  n/73) \big)\;.
\end{align*}
For $n>\frac{500}{\alpha}$, we have $\log n-\alpha  n/73<-2\log n$. Thus, $$\sum_{U  \subseteq V(H),  |U| \leq n /2, U\neq \emptyset } \Pr\left[E_{U}\right]<\sum_{k=1}^{\lceil\frac{n}2\rceil}  \exp(-2k \log n)<1\;,$$ as was needed.
\end{proof}
From now on, we fix an orientation of $\overrightarrow{H}$ as in Claim~\ref{cl:orbexist}.

Call a component $C$ of $F$ \emph{weakly untouched} if there is no oriented edge of $\overrightarrow{H}$ from $C$ to $V\setminus F$. Let $r'$ be the number of weakly untouched components. Obviously, each untouched component is also weakly untouched, and so $r\le r'$.
\begin{claim}
Let $C$ be an arbitrary component of $F$. Then $\Pr\left[\mbox{$C$ is weakly untouched}\right]\le \frac{1}{5}$.
\end{claim}
\begin{proof}
	Fix an $\alpha/3$-tough orientation of $H$. The probability that a component $C$ of $F$ is weakly untouched is at most
 \[ 
 \left(1-\frac{12}{ \alpha n}\right)^{\frac{\alpha}{3} |C| (n-|C|)}  \leq \frac{1}{e^{2}} \leq \frac{1}{5}\;. 
 \]
 \end{proof}
By the claim, we have $\Exp[r']\le cc(F)/5$. For each component of $F$ we have the event that that component is weakly untouched. Crucially, these events are independent.\footnote{This is not true for the collection of events of different components being untouched. The entire motivation for orienting the graph $F$ comes from this.} Hence, by Chernoff's bound,  $r^{'}\ge cc(F)/5+n^{4/5}$, with probability at most $n^{-4}$. It remains to argue that if this event $(*)$ does not occur, we have the outcome of the lemma. To this end, we observe that the number of components in the part of the graph that is formed by weakly touched components decreased by one half or even more (in the graph $F'$ compared to $F$). Hence, 
\[
cc(F')\le r'+0.5(cc(F)-r')=0.5(cc(F)+r')\leBy{$\neg$(*)} 0.6\cdot cc(F)+0.5n^{4/5},
\]
which finishes the proof.

\section{More on fractional multiples and $\kappa(W)$}\label{sec:proofsKappa}
\subsection{Bounds on and $\Upsilon(W)$ and $\kappa(W)$}
In this section, we establish several bounds about $\Upsilon(\cdot)$ and $\kappa(\cdot)$ used in the main proof.

First, we give a proof of Fact~\ref{fact:robustmindeg}, which together with Fact~\ref{fact:kappafinite} proves that $\kappa(W)<\infty$ for kernels $W$ arising in Theorem~\ref{thm:mainfull}.
\begin{proof}[Proof of Fact~\ref{fact:robustmindeg}]
We break the proof down into two parts. First, we prove that if $G$ is a $\rho$-robust graph, then for its graphon representation $W_G$ we have $\mindeg(W_G)\ge \rho/2$. Indeed, if $x\in\Omega$ is an element contained in a cell $\Omega_v\subset\Omega$ representing a vertex $v\in V(G)$, then the inequality $e_G(\{v\},V(G)\setminus \{v\})\ge\rho \cdot 1\cdot (n-1)\ge \rho n/2$ translates as $\deg_{W_G}(x)\ge \rho/2$.

Next, we prove that if $W$ is a cut distance limit of kernels $(W_n)_n$ with minimum degree at least $\rho/2$ then $\mindeg(W)\ge \rho/2$. This amounts to proving that for every $B\subset \Omega$ we have $\int_{B\times \Omega}W\ge \mu(B)\cdot \rho/2$. From the fact that $W$ is a limit of $(W_n)_n$ we have that 
\begin{equation}\label{eq:thankyou}
d_\square(W_n^{\pi_n},W)\to 0
\end{equation}
for a suitable sequence  $(\pi_n)_n$ of measure-preserving bijections. From the assumption on the minimum degrees of the kernels $W_n$ we have $\int_{B\times \Omega}W_n\ge \mu(B)\cdot \rho/2$. So, the claim follows from the definition of the cut norm distance, and from~\eqref{eq:thankyou}.
\end{proof}

We now prepare for the proof of Fact~\ref{fact:kappafinite}.
The next lemma deals with process $\fX_A$ and $\fX_B$ where the former dominates the latter. By \emph{domination} in this context we mean that the process can be coupled in a way that $\fX_A$ is a rooted supertree of $\fX_B$. While the extent of this definition may not be clear, we shall only use it in two instances. The first is when $A$ is a kernel with minimum degree at least $\beta$ and $B$ is the constant-$\beta$ kernel. The second is when $A=cB$ for $c\ge 1$.
\begin{lemma}\label{lem:UpsilonMonot} 
	Suppose that $A$ and $B$ are two kernels such that the process $\fX_A$ dominates the process $\fX_B$. Then we have $\Upsilon(A)\le \Upsilon(B)$.
\end{lemma}
\begin{proof}
	The lemma is implied immediately by the following trivial inequality: if $a_1,a_2,\ldots,b_1,b_2, \ldots \in[0,1]$ are sequences of numbers with $\sum_k a_k,\sum_k b_k\le 1$ and with $\sum_{k=1}^{L}a_k\le \sum_{k=1}^{L}b_k$ for each $L\in \NN$ then we have 
	\begin{equation*}
		\sum_{k=1}^\infty\frac{a_k}k\le \sum_{k=1}^\infty\frac{b_k}k\;.
	\end{equation*}
	Indeed, to finish the argument, we take $a_k:=\Pr[|\fX_A|=k]$ and $b_k:=\Pr[|\fX_B|=k]$.
\end{proof}

We now prove Fact~\ref{fact:kappafinite}. We break the proof down into two parts.
\begin{lemma}\label{lem:kappafinitedeg}
		\begin{enumerate}[label=(\roman*)]
		\item\label{en:kappaUpsilon} Suppose that $\beta\ge 1000$ and that $U$ is a kernel with minimum degree at least $\beta$. Then $\Upsilon(U)\le 5\beta^{-2}$.
		\item\label{en:kappakappa} Suppose that $\rho$ is positive and that $W$ is a kernel with minimum degree at least $\rho$. Then $\kappa(W)<\infty$.
	\end{enumerate}
\end{lemma}
\begin{proof}[Proof of~\ref{en:kappaUpsilon}]
The branching process $\fX_U$ dominates the Galton--Watson branching process with offspring distribution $\Poi(\beta)$. Lemma~\ref{lem:UpsilonMonot} tells us that $\Upsilon(U)\le \Upsilon(\beta)$.
So, it is enough to prove the bound for the constant-$\beta$ kernel. We need some large deviation bounds for the total progeny of Galton--Watson branching process with offspring distribution $\Poi(\beta)$. While there are plenty of bounds available in literature, a weak one suffices for our purposes, which we give a quick proof of. More specifically, we claim that we have
\begin{equation}\label{eq:GW100}
	\Pr\left[|\fX_\beta|<\frac{\beta^2}{4}\right]\le \beta^{-2}\;.
\end{equation}
To see that, first note that a short calculation (which we omit) gives that in this range the probability that a Poisson random variable with mean $\beta$ is less than $\frac\beta2$ is less than $\beta^{-3}$. We apply this fact first to the root particle of $\fX_\beta$ and see that typically it has at least $\frac\beta2$-many offspring. In that case we apply the same fact to the first $\lceil\frac\beta2\rceil$-many offspring and see that typically each of them has at least $\frac\beta2$-many offspring of its own. We conclude that, up to an error probability as on the right-hand side of~\eqref{eq:GW100}, the progeny of $\fX_\beta$ in the first three generations has already size at least $\frac{\beta^2}{4}$. Thus,
\begin{equation*}
\Upsilon(\beta)=\sum_{k=1}^\infty\frac{\Pr[|\fX_\beta|=k]}k\le \Pr\left[|\fX_\beta|<\frac{\beta^2}{4}\right]+\frac{\Pr\left[|\fX_\beta|\ge\frac{\beta^2}{4}\right]}{\frac{\beta^2}{4}}
\le \beta^{-2}+\frac1{\frac{\beta^2}{4}}=5\beta^{-2}\;.
\end{equation*}
\end{proof}
\begin{proof}[Proof of~\ref{en:kappakappa}]
We substitute into~\eqref{eq:defkappa},
\begin{align*}
	\kappa(W)&\le \int_{c=0}^{\frac{1000}{\rho}}1 \;\differential c+\int_{c=\frac{1000}{\rho}}^{+\infty}\sum_{k=1}^\infty \Upsilon(cW) \;\differential c\\
\JUSTIFY{by Part~\ref{en:kappaUpsilon}}	&\le \frac{1000}{\rho}+\int_{c=\frac{1000}{\rho}}^{+\infty} 5(c\rho)^{-2}\;\differential c<\infty\;.	
\end{align*}
\end{proof}

The next two lemmas tell us that the quantities $\Upsilon(\cdot)$ and $\kappa(\cdot)$ are almost the same for kernel $Z$ and for kernel $(1\pm \lambda)Z$ for small $\lambda$. The first lemma deals with $\Upsilon(\cdot)$.
\begin{lemma}\label{lem:MultiContinuityofUpsilon}
For every $\delta>0$, every $\lambda\in(0,\delta^2/16)$, for every kernel $Z$, we have $\Upsilon(Z)-\Upsilon((1-\lambda)Z)\le \delta$.
\end{lemma}
\begin{proof}
Let $L:=\lceil4/\delta\rceil$. For $k\in\NN$, we write $p_k:=\Pr[|\fX_{Z}|\le k]$ and $\tilde{p}_k:=\Pr[|\fX_{(1-\lambda)Z}|\le k]$. We have $\Upsilon(Z)=\sum_{k=1}^L \frac{p_k}{k^2}\pm \frac{1}{L}$ and $\Upsilon((1-\lambda )Z)=\sum_{k=1}^L \frac{\tilde{p}_k}{k^2}\pm \frac{1}{L}$.

Poisson thinning tells us that the process $\fX_{(1-\lambda)Z}$ can be obtained as follows. First, generate the tree $T$ corresponding to $\fX_{Z}$. Then we kill each nonroot vertex with probability $\lambda$ together with its entire progeny. The trimmed tree $T'$ we obtain in this way has distribution of $\fX_{(1-\lambda)Z}$. In particular, if we have the event that $v(T)\ge k$, then we can fix a subtree $F\subset T$ containing the root. If none of the nonroot vertices of $F$ got killed, then we have $v(T')\ge k$. This consideration gives $1-\tilde{p}_k\ge (1-\lambda)^{k-1}(1-p_k)$. Straightforward manipulations give
\begin{equation}\label{eq:MA}
\tilde{p}_k\le 1-(1-\lambda)^{k-1}(1-p_k)\le 1-(1-k\lambda)(1-p_k)\le k\lambda+p_k\;.
\end{equation}
We combine this with the expressions for $\Upsilon(Z)$ and $\Upsilon((1-\lambda)Z)$ above.
\begin{align*}
\Upsilon((1-\lambda )Z)-\Upsilon(Z)=\sum_{k=1}^{L}\frac{\tilde{p}_k-p_k}{k^2}\pm \frac{2}{L}\leByRef{eq:MA} \sum_{k=1}^{L}\frac{\lambda}{k}\pm \frac{2}{L}\le L\cdot \lambda\pm \frac{\delta}{2}<\delta\;,
\end{align*}
as was needed.
\end{proof}
\begin{lemma}\label{lem:MultiContinuityofKappa}
	For every $a>0$, and every kernel $W$ with positive minimum degree, there exists $b\in(0,1)$, so that we have $\kappa((1-b)W)\le \kappa(W)+a$ and $\kappa(W)\le \kappa((1+b)W)+a$.
\end{lemma}
\begin{proof}
Let $L\in(0,\infty)$ be such that $\int_{c=L}^\infty \Upsilon(\frac{c}2\cdot W)\differential c<a/4$. By Lemma~\ref{lem:kappafinitedeg}\ref{en:kappaUpsilon}, such an $L$ exists. By monotonicity (Lemma~\ref{lem:UpsilonMonot}), we have $\int_{c=L}^\infty \Upsilon(\alpha c W)\differential c<a/4$, for every $\alpha\ge \frac12$. Now, take $b:=a^2/(1600L^2)$. We have
\begin{align*}
	\kappa((1-b)W)&=\int_{c=0}^L\Upsilon((1-b)cW)\differential c\pm \frac{a}{4}\;,\\
	\kappa(W)&=\int_{c=0}^L\Upsilon(cW)\differential c\pm \frac{a}{4}\;,\\
	\kappa((1+b)W)&=\int_{c=0}^L\Upsilon((1+b)cW)\differential c\pm \frac{a}{4}\;.
\end{align*}
By Lemma~\ref{lem:MultiContinuityofUpsilon}, for the same $c$, the three integrands differ from each other by at most $a/(2L)$. Hence the lemma follows.
\end{proof}

\subsection{Computing $\kappa(\mathbbm{1})$}\label{ssec:ComputingZeta3}
In this section, we prove that $\kappa(\mathbbm{1})=\zeta(3)$. As we said, our argument is a version of an argument from~\cite{MR770868} recast to our language. Given the definition of $\Upsilon(\cdot)$, crucial to the proof is the total progeny distribution of the Galton--Watson process $\fX_c$, whose offspring distribution is $\Poi(c)$. The answer is given in the claim below.
\begin{claim*}
For every $c>0$ and each $k\in\NN$,  we have 
\[
\Pr[|\fX_c|=k]=c^{k-1}\cdot\exp(-ck)\frac{k^{k-1}}{k!}
\;.
\]
\end{claim*}
\begin{proof}
Let $\mathcal{P}_k$ be the set of isomorphism classes of all plane trees (recall Section~\ref{ssec:meaningofconstant}) of order $k$. We have $\Pr[|\fX_c|=k]=\sum_{P\in \mathcal{P}_k}\Pr[\fX_c\cong P]$. 

For each plane tree $P$ and for a vertex $v\in V(P)$, write $d_P(v)$ for the number of children of $v$, that is,
\begin{equation*}
	d_P(v)=\begin{cases}
		\deg_P(v)&\mbox{$v$ is the root of $P$,}\\
		\deg_P(v)-1&\mbox{otherwise.}
	\end{cases}
\end{equation*}
We write $\mathbf{d}_P$ for the collection $\{d_P(v)\}_{v\in V(P)}$. We use this notation in multinomial coefficients $\binom{k}{\mathbf{d}_P}$. For example, if $P\in \mathcal{P}_6$ is a plane tree in which the root has 3 children, of which the first and the third are leaves, and the second has 2 children (both of which are leaves), then we have $$\binom{6}{\mathbf{d}_P}=\binom{6}{3,0,2,0,0,0}\;.$$

The definition of the process $\fX_c$ gives that for each $P\in \mathcal{P}_k$, we have $$\Pr[\fX_c\cong P]=\prod_{v\in V(P)}\frac{c^{d_P(v)}\cdot\exp(-c)}{d_P(v)!}=\frac{c^{k-1}\cdot\exp(-ck)}{\prod_{v\in V(P)} d_P(v)!}\;.$$

For $P\in \mathcal{P}_k$, write $\mathcal{S}_P$ for the set of all pairs $(P,(L_v)_{v\in V(P)})$ where $(L_v)_{v\in V(P)}$ is a system of disjoint subsets of $[k]$ with the cardinality condition $|L_v|=d_P(v)$ for each $v\in V(P)$. This is like everyone in a multi-generational family buying cardigans for their own children (1 cardigan per 1 child, each cardigan is unique) but not being decided yet which child gets which cardigan. So, we have $|\mathcal{S}_P|=\binom{k}{\mathbf{d}_P}$. If $A\in \mathcal{S}_P$, then we write $V(A)$, $d_A(v)$, and $\mathbf{d}_A$ to mean $V(P)$, $d_P(v)$ and $\mathbf{d}_P$ respectively. Write $\mathcal{S}_k=\bigcup_{P\in\mathcal{P}_k}\mathcal{S}_P$. 

The above gives
\begin{align*}
\Pr[|\fX_c|=k]&=c^{k-1}\cdot\exp(-ck)\sum_{P\in \mathcal{P}_k}\frac{1}{\prod_{v\in V(P)} d_P(v)!}\\
&=c^{k-1}\cdot\exp(-ck)\sum_{A\in \mathcal{S}_k}\frac{1}{\binom{k}{\mathbf{d}_A}\prod_{v\in V(A)} d_A(v)!}
=c^{k-1}\cdot\exp(-ck)\sum_{A\in \mathcal{S}_k}\frac{1}{k!}\;.
\end{align*}
Last, observe that there is a natural bijection between each pair $(P,(L_v)_{v\in V(P)})\in\mathcal{S}_k$ and all rooted spanning trees on $K_k$. Indeed, the rooted spanning tree $T\subset K_k$ assigned to $(P,(L_v)_{v\in V(P)})\in\mathcal{S}_k$ is the unique tree isomorphic to $P$ whose vertices are as follows
\begin{itemize}
	\item the root of $T$ is the unique vertex not appearing in $\bigcup_{v\in V(P)}L_v$,
	\item each non-root vertex $v\in V(P)$ which is the $j$th child of its parent $w$ (with respect to the linear order on the vertices of $w$ in the plane tree $P$) is the $j$th smallest number in $L_w$.
\end{itemize}
By Cayley's formula there are $k^{k-1}$ rooted spanning trees of $K_k$. Hence,
\begin{align*}
	\Pr[|\fX_c|=k]=c^{k-1}\cdot\exp(-ck)\frac{k^{k-1}}{k!}\;,
\end{align*}
as was needed.
\end{proof}
We have
\begin{align*}
\kappa(1) &= \int_{c=0}^{+\infty}\Upsilon(c\cdot 1)\differential c=\int_{c=0}^{+\infty} \sum_{k=1}^{\infty} \frac{\Pr[|\fX_c|=k]}{k} \differential c\; \\
\JUSTIFY{by Claim}&=\int_{c=0}^{+\infty} \sum_{k=1}^{\infty}c^{k-1}\cdot\exp(-ck)\frac{k^{k-2}}{k!}\; \differential c=  \sum_{k=1}^{\infty}\int_{c=0}^{+\infty}c^{k-1}\cdot\exp(-ck)\frac{k^{k-2}}{k!}\; \differential c \\
\JUSTIFY{substitution $t=ck$}&=  \sum_{k=1}^{\infty}\int_{t=0}^{+\infty}\left(\frac{t}{k}\right)^{k-1}\cdot\exp(-t)\frac{k^{k-3}}{k!}\; \differential t
=  \sum_{k=0}^{\infty}\frac{1}{k^2\cdot k!}\int_{t=0}^{+\infty}t^{k-1}\cdot\exp(-t)\; \differential t \\
\JUSTIFY{Gamma function}&= \sum_{k=1}^{\infty}\frac{1}{k^2\cdot k!}\cdot (k-1)!=\zeta(3)\;.
\end{align*}

\begin{remark}
We believe that similar calculations should be tractable even for slightly more complicated kernels, such as $2\times 2$-step kernels, to obtain an explicit formula, but we have not pursued that direction.
\end{remark}

\section{Further questions}
\subsection{Structure of the minimum spanning tree}
In this paper we investigated the value of the random minimum spanning tree on dense graphs $(G_n)_n$ in which the distributions on the edges have well-behaved cumulative distribution functions with positive derivatives at~0. Our main theorem says that when the graphs $(G_n)_n$ with edges weighted by the derivatives of cumulative distribution function at~0 converge to a kernel $W$ then the values of the minimum spanning trees on $G_n$ converge in probability to a certain constant $\kappa(W)\in(0,\infty)$. This generalizes a classical result of Frieze about the random minimum spanning tree on complete graphs. The next natural step would be to study the geometry of the random minimum spanning tree. In the case of complete graphs, the corresponding results are due to Addario-Berry~\cite{LocalMST}, which determines the structure of the minimum spanning tree from a local perspective, and Addario-Berry--Broutin--Goldschmidt--Miermont~\cite{MR3706739} and Broutin--Marckert~\cite{ExplicitMST}, which determine the structure of the minimum spanning tree from a global perspective. Note that from both these perspectives, the minimum spanning trees behave substantially differently than the uniform spanning trees on complete graphs, whose local structure was determined in~\cite{Kolchin1977,Gri:RandomTree} and global in~\cite{MR1085326}. So, just like the results for the uniform spanning tree were recently extended to general sequences of dense graphs both in the local sense in~\cite{MR3876899} and in the global sense in~\cite{GlobalUSTGraphon}, we can ask for similar extensions of~\cite{LocalMST} and~\cite{MR3706739,ExplicitMST}.

\subsection{An extremal question for $\kappa(W)$}
What kernels $W$ have a small or a large $\kappa(W)$? In view of the fact that kernels with larger values tend to have smaller $\kappa$, it is natural to restrict to kernels of a fixed density. By linearity expressed by Fact~\ref{fact:AFM}, we can restrict to the case when the density is~1. Among robust graphons with $\|W\|_1=1$, we have $\sup_W \kappa(W)=\infty$. To see this, for each $\eps>0$ consider a kernel $W_\eps$ on $[0,1]^2$ defined by
\begin{align*}
	W_\eps(x,y)=\begin{cases}
		0 &\mbox{if $x,y\le\frac12$,}\\
		\eps &\mbox{if $x\le \frac12$ and $y>\frac12$\;,}\\
		4-2\eps &\mbox{if $x,y>\frac12$.}
	\end{cases}
\end{align*}
As $\eps\to0$, we have $\kappa(W_\eps)\to\infty$. While a calculation to show this is fairly straightforward, perhaps informal justification using random minimum spanning tree is more telling. So, suppose that we have a graph of order $n$ whose edges are equipped with probability distributions, so that the corresponding graph is close to $W_\eps$. For the $\approx\frac{n}2$ many vertices $v$ corresponding to~$[0,\frac12]$, the minimum \previouslyweight\ of an edge incident to $v$ is $\approx \frac2{\eps n}$ in expectation. While this value is not concentrated at a particular vertex, we can use the law of large numbers to deduce concentration of the sum. That is, incorporating these vertices into a spanning tree costs $\approx \frac{n}2\cdot\frac2{\eps n}=\frac{1}{\eps}$. Hence indeed, $\sup_{W:\|W\|_1=1} \kappa(W)=\infty$.

So, the meaningful question is about a lower bound. We believe that the optimal bound comes from 1-regular kernels.
\begin{conjecture}
Suppose that $W$ is a kernel with $\|W\|_1=1$. Then $\kappa(W)\ge \zeta(3)$.
\end{conjecture} 

\subsection{Other optimization problems with random inputs}
There are many other optimization problems on randomly weighted graphs. Probably the most famous and closest to our random minimum spanning tree problem is the random assignment problem: \emph{In a complete graph $K_n$ ($n$ is even), whose edges are equipped with $\Uni[0,1]$-\previouslyweights, find the total length of a minimum perfect matching.}\footnote{This is actually called `minimum weight perfect matching' in literature, the only reason we use the word `length' is to stay consistent with Footnote~\ref{foot:weightlenght}.} Analysing a minimum length total perfect matching is a substantially more challenging problem, in particular because no greedy approach such as the one of Kruskal's algorithm exists. Rigorously establishing a physics prediction by Mézard and Parisi~\cite{MezardParisi}, Aldous~\cite{MR1839499} proved that the weight of a minimum perfect matching converges in probability to $\zeta(2)/2$. Another optimization problems involve Hamiltonian cycles and 2-factors, see~\cite{MR2104159,MR2600434}. All these optimization problems call for inhomogeneous counterparts.

\section*{Acknowledgments}
The problem considered in this paper stems from discussions with Ellie Archer and Matan Shalev. Ellie and Matan also provided helpful comments to an early version of this manuscript. We thank two referees for their comments.

\bibliographystyle{plain}
\bibliography{DSG}

\begin{thebibliography}{10}

\bibitem{LocalMST}
Louigi Addario-Berry.
\newblock The local weak limit of the minimum spanning tree of the complete
  graph.
\newblock arXiv:1301.1667.

\bibitem{MR3706739}
Louigi Addario-Berry, Nicolas Broutin, Christina Goldschmidt, and Gr\'{e}gory
  Miermont.
\newblock The scaling limit of the minimum spanning tree of the complete graph.
\newblock {\em Ann. Probab.}, 45(5):3075--3144, 2017.

\bibitem{MR1085326}
David Aldous.
\newblock The continuum random tree. {I}.
\newblock {\em Ann. Probab.}, 19(1):1--28, 1991.

\bibitem{MR2023650}
David Aldous and J.~Michael Steele.
\newblock The objective method: probabilistic combinatorial optimization and
  local weak convergence.
\newblock In {\em Probability on discrete structures}, volume 110 of {\em
  Encyclopaedia Math. Sci.}, pages 1--72. Springer, Berlin, 2004.

\bibitem{MR1839499}
David~J. Aldous.
\newblock The {$\zeta(2)$} limit in the random assignment problem.
\newblock {\em Random Structures Algorithms}, 18(4):381--418, 2001.

\bibitem{GlobalUSTGraphon}
Eleanor Archer and Matan Shalev.
\newblock The {GHP} scaling limit of uniform spanning trees of dense graphs.
\newblock {\em Random Structures Algorithms}, 65(1):149--190, 2024.

\bibitem{MR1721947}
Andrew Beveridge, Alan Frieze, and Colin McDiarmid.
\newblock Random minimum length spanning trees in regular graphs.
\newblock {\em Combinatorica}, 18(3):311--333, 1998.

\bibitem{MR2599196}
Béla Bollob\'{a}s, Christian Borgs, Jennifer Chayes, and Oliver Riordan.
\newblock Percolation on dense graph sequences.
\newblock {\em Ann. Probab.}, 38(1):150--183, 2010.

\bibitem{MR2337396}
Béla Bollob\'{a}s, Svante Janson, and Oliver Riordan.
\newblock The phase transition in inhomogeneous random graphs.
\newblock {\em Random Structures Algorithms}, 31(1):3--122, 2007.

\bibitem{MR2455626}
C.~Borgs, J.~T. Chayes, L.~Lov\'{a}sz, V.~T. S\'{o}s, and K.~Vesztergombi.
\newblock Convergent sequences of dense graphs. {I}. {S}ubgraph frequencies,
  metric properties and testing.
\newblock {\em Adv. Math.}, 219(6):1801--1851, 2008.

\bibitem{ExplicitMST}
Nicolas Broutin and Jean-François Marckert.
\newblock Convex minorant trees associated with {B}rownian paths and the
  continuum limit of the minimum spanning tree.
\newblock arXiv:2307.12260.

\bibitem{DevFra:Connectivity}
Luc Devroye and Nicolas Fraiman.
\newblock Connectivity of inhomogeneous random graphs.
\newblock {\em Random Structures Algorithms}, 45(3):408--420, 2014.

\bibitem{MR770868}
A.~M. Frieze.
\newblock On the value of a random minimum spanning tree problem.
\newblock {\em Discrete Appl. Math.}, 10(1):47--56, 1985.

\bibitem{MR1054012}
A.~M. Frieze and C.~J.~H. McDiarmid.
\newblock On random minimum length spanning trees.
\newblock {\em Combinatorica}, 9(4):363--374, 1989.

\bibitem{MR2104159}
Alan Frieze.
\newblock On random symmetric travelling salesman problems.
\newblock {\em Math. Oper. Res.}, 29(4):878--890, 2004.

\bibitem{MR4482093}
Jan Greb\'{i}k and Israel Rocha.
\newblock Fractional isomorphism of graphons.
\newblock {\em Combinatorica}, 42(3):365--404, 2022.

\bibitem{Gri:RandomTree}
G.~R. Grimmett.
\newblock Random labelled trees and their branching networks.
\newblock {\em J. Austral. Math. Soc. Ser. A}, 30(2):229--237, 1980/81.

\bibitem{MR3876899}
Jan Hladk\'{y}, Asaf Nachmias, and Tuan Tran.
\newblock The local limit of the uniform spanning tree on dense graphs.
\newblock {\em J. Stat. Phys.}, 173(3-4):502--545, 2018.

\bibitem{HlaHngLim:GraphonBranching}
Jan Hladký, Eng~Keat Hng, and Anna~Margarethe Limbach.
\newblock Graphon branching processes and fractional isomorphism.
\newblock arXiv:2408.02528.

\bibitem{MR2816939}
Svante Janson and Oliver Riordan.
\newblock Duality in inhomogeneous random graphs, and the cut metric.
\newblock {\em Random Structures Algorithms}, 39(3):399--411, 2011.

\bibitem{MR2880662}
Svante Janson and Oliver Riordan.
\newblock Susceptibility in inhomogeneous random graphs.
\newblock {\em Electron. J. Combin.}, 19(1):Paper 31, 59, 2012.

\bibitem{Kolchin1977}
V.~F. Kolchin.
\newblock Branching processes, random trees, and a generalized scheme of
  arrangements of particles.
\newblock {\em Mathematical notes of the Academy of Sciences of the USSR},
  21(5):386--394, May 1977.

\bibitem{KruskalAlgorithms}
Joseph~B. Kruskal, Jr.
\newblock On the shortest spanning subtree of a graph and the traveling
  salesman problem.
\newblock {\em Proc. Amer. Math. Soc.}, 7:48--50, 1956.

\bibitem{Lovasz2006}
László Lov{\'a}sz and Balázs Szegedy.
\newblock {Limits of dense graph sequences}.
\newblock {\em J. Combin. Theory Ser. B}, 96(6):933--957, 2006.

\bibitem{Lovasz2012}
László Lovász.
\newblock {\em {Large networks and graph limits}}, volume~60 of {\em {American
  Mathematical Society Colloquium Publications}}.
\newblock American Mathematical Society, Providence, RI, 2012.

\bibitem{MezardParisi}
M.~Mézard and G.~Parisi.
\newblock On the solution of the random link matching problem.
\newblock {\em J. Physique}, 48(9):1451--1459, 1987.

\bibitem{MR0905183}
J.~Michael Steele.
\newblock On {F}rieze's {$\zeta(3)$} limit for lengths of minimal spanning
  trees.
\newblock {\em Discrete Appl. Math.}, 18(1):99--103, 1987.

\bibitem{MR2600434}
Johan W\"{a}stlund.
\newblock The mean field traveling salesman and related problems.
\newblock {\em Acta Math.}, 204(1):91--150, 2010.

\end{thebibliography}

\end{document}